\documentclass[11pt,a4paper]{amsart}

\usepackage{amsmath,amsthm}
\usepackage{amsfonts}
\usepackage{latexsym}
\usepackage{amsbsy}
\usepackage{amssymb}
\usepackage{verbatim}
\usepackage{graphicx}
\usepackage{hyperref}

\def\NN{\mathbb{N}}

\def\RR{\mathbb{R}}
\def\CC{\mathbb{C}}
\def\PP{\mathbb{P}}

\def\EE{\mathbf{E}}
\def\P{\mathbf{P}}
\def\L{\mathcal{L}}

\def\D{\mathcal{D}}
\def\F{\mathcal{F}}
\def\G{{\mathcal G}}

\def\1{\mathbf{1}}
\def\0{\mathbf{0}}

\def\I{\mathbf{I}}

\newcommand{\te}{\widehat{e}}

\newtheorem{Thm}{Theorem}

\newtheorem{Prop}{Proposition}

\def\BEN{\begin{enumerate}}  \def\BI{\begin{itemize}}
\def\EEN{\end{enumerate}}   \def\EI{\end{itemize}}

\def\mc{\mathcal}  \def\ovl{\overline}

\def\le{\left}
\def\ri{\right}
\def\te#1{\mathrm{e}^{#1}}   
\def\WT{\widetilde}
\def\WH{\widehat} 
   \def\ii{\mathrm{i}}

\def\F{\Phi} 
  \def\td{\text{\rm d}}
 
\begin{document}
\title{On additive time-changes of Feller processes}
\author{Aleksandar Mijatovi\'{c}}
\address{Department of Mathematics, Imperial College London}
\email{a.mijatovic@imperial.ac.uk}
\author{
Martijn Pistorius}
\address{Department of Mathematics, Imperial College London}
\email{m.pistorius@imperial.ac.uk}

\begin{abstract}
In this note we generalise the Phillips theorem~\cite{Phillips}
on the subordination of Feller processes by L\'evy subordinators
to the class of additive subordinators (i.e. subordinators with
independent but possibly nonstationary increments).
In the case where the original Feller process is L\'evy
we also express the time-dependent characteristics of 
the subordinated process 
in terms of the characteristics of the L\'evy process and
the additive subordinator.
\end{abstract}


\maketitle

\section{Introduction} 
One of the established devices  
for building statistically relevant market models is that of the stochastic  
change of time-scale (e.g. Carr et al.~\cite{CarrMadanGemanYor_selfdec}). 
Such a time change may be modelled as an independent additive
subordinator 
$Z=\{Z_t\}_{t\geq0}$, i.e.
an increasing stochastic process with independent
possibly nonstationary
increments. 
If we subordinate a time-homogeneous Markov process 
$X=\{X_t\}_{t\geq0}$
by 
$Z$, 
the resulting process 
$Y=\{X_{Z_t}\}_{t\geq0}$
is a Markov process that will in general be 
time-inhomogeneous. 
The main result of this note shows that 
if
$X$
is a Feller process and
$Z$
satisfies some regularity assumptions, 
then
$Y$
is a time-inhomogeneous 
Feller process. 
The generator of 
$Y$
is expressed 
in terms of the generator of 
$X$
and the characteristics of 
$Z$.
In the special case where 
$X$
is a L\'evy process 
it is shown that
$Y$
is an additive process with characteristics 
that 
are given 
explicitly in terms of the characteristics of 
$X$
and of the additive subordinator
$Z$.
The explicit knowledge of the 
generator 
of 
$Y$
is desirable from the viewpoint of pricing theory
because contingent claims 
in the time-inhomogeneous
market model
$Y$
can be evaluated using algorithms 
that are based on the explicit form of the generator of the 
underlying process (see for example \cite{MG_Alg}).

\section{Time-changed Feller processes} 

Throughout the paper we assume that
$X=\{X_t\}_{t\geq0}$ 
is a c\`adl\`ag
Feller process with the state-space
$\RR^n$
for some
$n\in\NN$
and the infinitesimal generator 
$\L$
defined on a dense subspace
$\D(\L)$
in the Banach space of all continuous functions
$C_0(\RR^n)$
that vanish at infinity with norm
$\|f\|_\infty:=\sup_{x\in\RR^n}|f(x)|$.
The corresponding semigroup 
$(P_t)_{t\ge 0}$
is given  by
$P_tf(x)=\EE^x[f(X_t)]$
for any
$f\in C_0(\RR^n)$,
where the expectation is taken with respect to 
the law of 
$X$
started at
$X_0=x$
(see Ethier and Kurtz~\cite{EthierKurtz} for the definition
and properties of
Feller semigroups).

Let
$Z=\{Z_t\}_{t\geq0}$ 
be an additive process,
independent of
$X$,
with the Laplace exponent
$\ovl\psi_t(u) = \log\EE[\te{-u Z_t}]$
given by $\ovl\psi_t(u):= \int_0^t \psi_s(u)\td s$, where
$\beta:\RR_+\to\RR_+$,
$g:\RR_+\times\mathbb R_+\to\RR_+$
are
continuous 
and for all 
$s\in\RR_+$,
$u\in\CC$
we have
$\int_{(0,\infty)}(1\wedge r)g(s,r)\td r<\infty$
and 
\begin{eqnarray}
\label{eq:psi}
\psi_s(u)  & = &  -u\beta(s) + \int_{(0,\infty)}(\te{-ur} - 1)g(s,r)\td r\quad\text{if}\quad
\Re(u)\geq0.
\end{eqnarray}
In other words 
$Z$
is a c\`adl\`ag process with nondecreasing paths such that
the random variable
$Z_t-Z_s$
is independent of 
$Z_u$
for all
$0\leq u\leq s<t$
(see
Jacod and Shiryaev~\cite{JacodShiryaev} Ch. II, Sec. 4c,
for a systematic treatment of additive processes).

In this paper we are interested in the process
$(D,Y)=\{(D_t,Y_t)\}_{t\geq 0}$
defined by
$D_t:=D_0+t$ and $Y_t := X_{Z_{D_t}}$ 
for some
$D_0\in\RR_+$.

\begin{Thm} 
\label{thm:TimeDeptPhillips}
The process 
$(D,Y)$ 
is Feller 
with the state-space~$\RR_+\times\RR^n$ 
and
infinitesimal generator $\L'$,
defined on a dense subspace of 
the Banach space
$C_0(\RR_+\times\RR^n)$
of continuous functions that vanish at infinity,
given~by
$$
\L'f(s,x) = \frac{\partial{f}}{\partial s}(s,x) +
\beta(s)\L f_s(x) + \int_{(0,\infty)}[P_r f_s(x) - f(s,x)]\,g(s,r)\td r,
$$
where 
$f\in C_0(\RR_+\times\RR^n)$
such that 
$f_s(\cdot):=f(s,\cdot)\in \D(\L)$
$\forall s\in\RR_+$
and the functions
$(s,x)\mapsto \mathcal Lf_s(x)$ and $(s,x)\mapsto \frac{\partial{f}}{\partial s}(s,x)$
are continuous and vanish at infinity.
\end{Thm}
If 
$Z$ 
is a L\'{e}vy subordinator, Theorem~\ref{thm:TimeDeptPhillips} 
reduces to the well-known Philips~\cite{Phillips} theorem. 
If
$X$ 
is a L\'{e}vy process, then the time-changed 
process is an additive process with
characteristics determined by those of 
$Z$
and 
$X$.

\begin{Prop}
\label{cor:timechange}
Let $X$ be a L\'{e}vy process 
with $X_0=0$
and 
characteristic triplet 
$(c,Q, \nu)$,  
where
$c\in\RR^n$,
$Q\in\RR^{n\times n}$
a nonnegative symmetric matrix and
$\nu$
a measure on
$\RR^n\backslash \{0\}$
such 
that 
$\int_{\RR^n\backslash \{0\}} (|x|^2\wedge 1)\nu(\td x)$.
The process
$Y$
defined above (with
$D_0=0$)
is additive with c\`adl\`ag paths,
jump measure 
$$
\WT\nu_s(\td x) =\beta(s)\nu(\td x) +
\int_{(0,\infty)}\P(X_r\!\in\!\td x)g(s,r)\td r,$$
nonnegative symmetric matrix 
$ \WT Q_s =  \beta(s)Q$, 
drift 
$$\WT c_s = \beta(s)c + \int_{(0,\infty)}\EE[X_r\I_{\{|X_r|\leq 1\}}]g(s,r)\td r$$
and characteristic exponent
$\ovl\Psi_t(u)=\int_0^t\Psi_s(u) \td s$
(recall that
$\EE[\te{\ii u\cdot Y_t}]=\te{\ovl\Psi_t(u)}$
for all
$u\in\RR^n$)
where
$$\Psi_s(u) = 
\ii u\cdot \WT c_s -\frac{1}{2}u\cdot \WT Q_s u +\int_{\RR^n\backslash \{0\}}\left[\te{\ii u\cdot x}-1-\ii (u\cdot x) \I_{\{|x|\leq1\}}\right]\WT\nu_s(\td x).$$
\end{Prop}

\section{Example: a symmetric self-decomposable process}

Suppose that $Y$ is an additive process, considered in~\cite{CarrMadanGemanYor_selfdec}
as a model for the risky security,
with no drift or Gaussian component and jump density
$$
g_Y(t,y) = h_\nu(|y|/t^\gamma) \frac{\gamma}{\nu t^{\gamma + 1}},\quad\text{where}\quad
h_\nu(y) = \frac{1}{\nu} \exp(-y/\nu)\I_{\{y>0\}}.
$$
Then in law  the process
$Y$ 
is equal to a Brownian motion time-changed by 
an independent additive subordinator 
$Z$ 
with 
$\beta\equiv0$
and jump density
$$
g(t,r) 
= a_t\te{-r/b_t},\quad\text{where}\quad
a_t=\frac{\gamma}{\nu^3t^{2\gamma + 1}},\quad
b_t=2\nu^2t^{2\gamma}.
$$
It is clear from Proposition~\ref{cor:timechange} 
that
$\WT c_t=\WT Q_t=0$
for all
$t\in\RR_+$
and that
the moment-generating functions of measures
$\WT\nu_t(\td x)$
and 
$g_Y(t,x)\td x$
coincide
\begin{equation*}
\int_{\RR\backslash\{0\}} \te{\lambda x}\>\WT\nu_t(\td x)=
\frac{2\gamma}{\nu t(1 - \lambda^2 \nu^2t^{2\gamma})}=
\int_{\RR\backslash\{0\}} \te{\lambda x}g_Y(t,x)\td x  
\end{equation*}
for  
$|\lambda|<1/\nu t^\gamma$.
This implies that the two additive processes coincide in~law.

\section{Proofs}
\subsection{Proof of Proposition~\ref{cor:timechange}}
Let
$\Psi_X(u)$
denote the characteristic exponent of 
the L\'evy process
$X$,
i.e.
$ \EE[\exp(\ii u\cdot X_s)]=\exp(s\Psi_X(u))$
for any
$u\in\RR^n$.
Since
$X$ and $Z$
are independent processes with independnent increments,
for any sequence of positive real numbers
$0\leq t_0< \ldots < t_m$ 
and 
vectors
$u_1, \ldots, u_m\in\RR^n$
it follows that
\begin{eqnarray*}
\EE\left[\te{\ii \sum_{i=1}^m u_i\cdot (Y_{t_i} -
Y_{t_{i-1}})}\right] &=& \EE\left[\EE\Big[\prod_{i=1}^m \te{\ii u_i\cdot 
(X_{Z_{t_i}} -
X_{Z_{t_{i-1}}})}\Big\arrowvert Z_{t_0},\ldots,Z_{t_m} \Big]\right]\\
&=& \EE\left[\prod_{i=1}^m \te{ (Z_{t_i}- Z_{t_{i-1}})\Psi_X(u_i)}\right] \\
& = & \prod_{i=1}^m \EE\left[\te{\ii u_i\cdot 
(Y_{t_i} - Y_{t_{i-1}})}\right].
\end{eqnarray*}
Hence the process
$Y$ 
also has independent increments. 
Since 
$Y$
is clearly c\`adl\`ag
(as $X$ and $Z$ are),
it is an additive process.

Finally, we have to determine the characteristic curve of $Y$.
An argument similar to the one above implies that the characteristic function of 
$Y_t$ 
equals
\begin{eqnarray*}
\EE[\te{\ii  u\cdot Y_t}] &=&  \EE[\te{\Psi_X(u)Z_t}]=
\te{\int_0^t \psi_s(-\Psi_X(u))\td s} \quad\text{for any}\quad u\in\RR^n.
\end{eqnarray*}
The last equality holds since
$\Re(\Psi_X(u))\leq 0$
for all
$u$
and the integral in~\eqref{eq:psi}
is well-defined.
It is not difficult to prove that
for any L\'evy process
$X$
started at
$0$
there exists a constant
$C>0$
such that the inequality holds
$$\max\left\{\P(|X_r|>1), |\EE[X_r\I_{\{|X_r|\leq 1\}}]|,\EE[|X_r|^2\I_{\{|X_r|\leq 1\}}]\right\}\leq C (r\wedge1)\quad\forall r\in\RR_+$$
(see e.g. Lemma 30.3 in  Sato~\cite{Sato}). Therefore, since
$\int_0^\infty\!\!g(s,r)(r\wedge 1)\td r<\infty$
by assumption,
we have
$$\int_0^\infty\! g(s,r)\max\left\{\P(|X_r|>1), |\EE[X_r\I_{\{|X_r|\leq 1\}}]|,\EE[|X_r|^2\I_{\{|X_r|\leq 1\}}]\right\}
\td r<\infty.$$
We can thus define the measure 
$\WT\nu_s(\td x)$,
the vector
$\WT c_s$,
the matrix
$\WT Q_s$
and
the function
$\Psi_s(u)$
by the formulae in 
Proposition~\ref{cor:timechange}.
The L\'evy-Khintchine representation
$$\Psi_X(u) = 
\ii u\cdot c -\frac{1}{2}u\cdot Q u +\int_{\RR^n\backslash \{0\}}\left[\te{\ii u\cdot x}-1-\ii (u\cdot x) \I_{\{|x|\leq1\}}\right]\nu(\td x)$$
and Fubini's theorem, which applies by the inequality above, yield the following
calculation, whcih concludes the proof of the proposition:
\begin{eqnarray*}
\psi_s(-\Psi_X(u)) &=& \beta(s) \Psi_X(u) + 
\int_0^\infty\!\!\!(\EE[\te{\ii u\cdot X_r}] - 1)g(s,r)\td r\\
& = &  \beta(s) \Psi_X(u) + \ii u\cdot \int_0^\infty\!\!\!\EE[X_r\I_{\{|X_r|\leq 1\}}]g(s,r)\td r  \\
& +  & 
\int_0^\infty\!\!\!(\EE[\te{\ii u\cdot X_r}] - 1 - \ii u\cdot \EE[X_r\I_{\{|X_r|\leq 1\}}])g(s,r)\td r
 =  \Psi_s(u). 
%
\end{eqnarray*}

\subsection{Proof of Theorem~\ref{thm:TimeDeptPhillips}}
Note first that the paths of the process
$(D,Y)$
are c\`adl\`ag.
In what follows
we prove that 
$(D,Y)$
is a Markov process that satisfies the Feller property
and find the generator of its semigroup.

\noindent {\bf 1. Markov property.} 
For any 
$g\in C_0(\RR_+\times \RR^n)$ 
define
\begin{eqnarray*}
Q_tg(s,x) &:= & \EE[g(D_t, Y_t)|D_0=s,Y_0=x] =
\EE[g(s+t,X_{Z_{s+t}})|X_{Z_s}=x].
\end{eqnarray*}
Let
$\lambda_{s,s+t}(\td r) := \PP(Z_{s+t}-Z_s\in\td r)$ 
denote the law of the increment of 
$Z$
which may have an atom at $0$. 
Then, since
$X$
and
$Z$
are independent processes and
the increments of 
$Z$
are independent of the past, it follows from the definition that
$$Q_tg(s,x) = \int_{[0,\infty)} \EE[g(s+t,X_{Z_{s}+r})|X_{Z_s}=x]\lambda_{s,s+t}(\td r).$$
Define for any
$v\in\RR_+$
a
$\sigma$-algebra
$\G_v=\sigma(X_l: l\in[0,v])$.
Then 
for a Borel set $A\in\mc B(\RR^n)$ and any 
$X_0=x_0\in\RR^n$ 
the Markov property of 
$X$ 
yields
\begin{eqnarray*}
\EE^{x_0}[g(t+s, X_{Z_s+r})\I_{\{X_{Z_s}\in A\}}] &\!\!=& \!\!\!
\int_{\![0,\infty)}\!\!\!\!\!\!\!\! \EE^{x_0}[g(t+s, X_{v+r})\I_{\{X_v\in A\}}] \lambda_{0,s}(\td v)\\
&\!\!=& \!\!\!
\int_{\![0,\infty)}\!\!\!\!\!\!\!\! \EE^{x_0}\left[\EE[g(t+s, X_{v+r})|\mc G_v]\I_{\{X_v\in A\}}\right]\!\!\lambda_{0,s}(\td v)\\
&\!\!=&
\EE^{x_0}\left[\EE^{X_{Z_s}}[g(t+s,X_r)]\I_{\{X_{Z_s}\in A\}}\right].
\end{eqnarray*}
Hence we get 
$\EE[g(t+s,X_{Z_s+r})|X_{Z_s}] = \EE^{X_{Z_s}}[g(t+s,X_r)]$~a.s. 
for any $r\in\RR_+$ 
and the following identity holds
\begin{equation}\label{eq:Qt}
Q_tg(s,x) = \int_{[0,\infty)} \EE^x[g(s+t,X_r)]\lambda_{s,s+t}(\td r).
\end{equation}
A similar argument and the monotone class theorem 
imply that,
if
$\F_s=\sigma(X_{Z_l}: l\in[0,s])$, 
then
$$\EE[g(t+s,X_{r+Z_s})|\F_s] = \EE^{X_{Z_s}}[g(t+s,X_r)]\quad\text{a.s.}$$
The process
$(D,Y)$,
started at
$(0,x_0)$,
satisfies
$$\EE[g(D_{s+t}, Y_{s+t})|\mc F_s] = \EE[g(t+s,X_{Z_{s+t}})|\mc F_s] = 
Q_tg(s,X_{Z_s}) = Q_tg(D_s,Y_s)$$
and is 
therefore
Markov 
with the semigroup
$(Q_t)_{t\geq0}$.

\noindent {\bf 2. Feller property.} 
Since
$(D,Y)$
and
$Z$
are right-continuous,
identity~\eqref{eq:Qt}
implies that 
$\lim_{t\searrow0}Q_tf(s,x)=f(s,x)$ 
for each 
$(s,x)\in\RR_+\times\RR^n$.  
It is well-known that in this case pointwise convergence implies
convergence in the Banach space
$(C_0(\RR_+\times \RR^n),\|\cdot\|_\infty)$.
It also follows from 
representation~\eqref{eq:Qt},
the dominated convergence theorem and the
Feller property of $X$ 
that a continuous function
$(s,x)\mapsto Q_tg(s,x)$ 
tends to zero at infinity
for any
$g\in C_0(\RR_+\times \RR^n)$.
Hence $(D,Y)$ is a Feller process.

\noindent {\bf 3. Infinitesimal generator of the semigroup $(Q_t)_{t\geq0}$.}
As before,
let
$\lambda_{s,s+t}$
be the law of the increment 
$Z_{s+t}-Z_s$
and let
$\psi_s$
be as in~\eqref{eq:psi}. 
Let
$(t_n)_{n\in\NN}$
be a decreasing sequence in
$(0,\infty)$
that converges to zero.
Denote 
by
$\WH\mu_n$
the Laplace transform of a compound Poisson process
with L\'evy measure
$t_n^{-1}\lambda_{s,s+t_n}$.
Hence we find for any
$u\in\CC$
that satisfies
$\Re(u)\geq0$
\begin{eqnarray*}
\WH\mu_n(u) &=& \exp\le(\frac{1}{t_n}\int_0^\infty (\te{-u r} - 1) \lambda_{s,s+t_n}(\td r)\ri)\\
& =  & \exp\le(t_n^{-1}(\te{\int_s^{s+t_n}\psi_v(u)\td v} - 1)\ri).
\end{eqnarray*}
Since the function
$t\mapsto \int_0^t \psi_{s+v}(u)\td v$ 
is right-differentiable at zero with derivative $\psi_s(u)$, 
we get
\begin{equation*}
\lim_{n\to\infty}\WH\mu_n(u) = \exp(\psi_s(u)).
\end{equation*}
It is clear from~\eqref{eq:psi}
that
$\exp(\psi_s(u))$ 
is a Laplace transform of an infinitely divisible distribution
with L\'evy measure
$g(s,r)\td r$.
Therefore by Theorem~8.7 in~\cite{Sato} 
for every continuous bounded function 
$k:\RR\to\RR$
that vanishes on a neighbourhood of zero
we get 
\begin{equation}\label{eq:conv1}
\lim_{n\to\infty}t_n^{-1}\int_0^\infty k(r) \lambda_{s,s+t_n}(\td r) = \int_0^\infty k(r)g(s,r)\td r.
\end{equation}
Furthermore the same theorem implies that 
for any continuous function $h$ such that $h(r)=1+o(|r|)$ for 
$|r|\to 0$ 
and 
$h(r)=O(1/|r|)$ 
for 
$|r|\to\infty$
we have
\begin{equation}\label{eq:conv2}
\lim_{n\to\infty}t_n^{-1}\int_0^\infty rh(r) \lambda_{s,s+t_n}(\td r) = \beta(s) + 
\int_0^\infty rh(r)g(s,r)\td r.
\end{equation}
A key observation is that~\eqref{eq:conv1}
and~\eqref{eq:conv2}
together imply that~\eqref{eq:conv1}
holds for every continuous 
bounded function 
$k$
that satisfies
$k(r)=o(|r|)$
as
$r\searrow0$.

\noindent \textbf{Claim.} Let the function 
$f\in C_0(\RR_+\times\RR^n)$ 
satisfy the assumptions of Theorem~\ref{thm:TimeDeptPhillips}.
Then for any 
$(s,x)\in\RR_+\times\RR^n$
the limit holds
\begin{eqnarray*}
\lim_{t\searrow0}t^{-1} (Q_t f - f)(s,x) &=&
 \frac{\partial f}{\partial s}(s,x) + \beta(s) \L f_s(x)\\ 
 &+& \int_0^\infty\!\!\![P_rf_s(x) - f_s(x)]g(s,r)\td r.
\end{eqnarray*}

To prove this claim
recall first that
$(P_t)_{t\geq0}$
is the semigroup of 
$X$
and note that the identity holds
$$ (Q_t f - f)(s,x) =
\EE^{s,x}[f(D_t,Y_t) - f(D_0,Y_t)] + \int_0^\infty[P_rf_s(x) - f_s(x)]\lambda_{s,s+t}(\td r).$$
If we divide this expression by 
$t$
and take the limit as
$t\searrow0$,
the first term converges to the partial derivative 
$\frac{\partial f}{\partial s}(s,x)$ 
by the dominated convergence theorem
(recall that the paths of 
$Y$
are right-continuous). 

Choose a function
$h$ 
as above,
define
$D(r):=P_rf_s(x) - f_s(x)$
and express the second term as
\begin{eqnarray*}
t^{-1}\int_0^\infty D(r) \lambda_{s,s+t}(\td r) & = & 
\lefteqn{t^{-1}\int_0^\infty D(r)(1-h(r)) \lambda_{s,s+t}(\td r)}\\
&+& t^{-1}\int_0^\infty (D(r) - r \L f_s(x)) h(r)\lambda_{s,s+t}(\td r) \\
&+ & \L f_s(x) t^{-1} \int_0^\infty rh(r)\lambda_{s,s+t}(\td r).
\end{eqnarray*}
The first and second integrals on the right-hand side 
converge by~\eqref{eq:conv1} 
to
$$
\int_0^\infty D(r)(1-h(r)) g(s,r)\td r\quad\text{and}\quad
\int_0^\infty (D(r) - r \L f_s(x)) h(r)g(s,r)(\td r) 
$$
respectively and the third integral converges 
by~\eqref{eq:conv2} 
to
$$
\L f_s(x)  \int_0^\infty rh(r)g(s,r)\td r.
$$
This proves the claim.

Since 
$(Q_t)_{t\geq0}$
is a strongly continuous contraction semigroup
on the function space
$C_0(\RR_+\times\RR^n)$
with some generator
$\L'$,
if the pointwise limit in the claim exists
and is in
$C_0(\RR_+\times\RR^n)$
for some  continuous function 
$f$
that vanishes at infinity,
then 
$f$
is in the domain of 
$\L'$
and 
$\L'f$
equals this limit
(see e.g. Lemma~31.7 in \cite{Sato}). 
This concludes the proof of the theorem.

\end{document}